\documentclass{article}
\usepackage{a4wide,hyperref}

\usepackage{graphicx}
\usepackage{amsmath, amssymb}
\usepackage{booktabs}
\usepackage{pstricks}
\usepackage{pst-node}
\usepackage{xspace}

\newcommand{\cip}{\mbox{$\,\perp\!\!\!\perp\,$}}

\newcommand{\ie}{{\em i.e.\/}\xspace}

\newcommand{\secref}[1]{\S$\,$\ref{sec:#1}}
\newcommand{\tabref}[1]{Table~\ref{tab:#1}}
\newcommand{\figref}[1]{Figure~\ref{fig:#1}}

 \newcommand{\eg}{{\em
    e.g.\/}\xspace} 
\newcommand{\pc}{\mbox{\rm PC}}

\renewcommand{\b}[1]{\textbf{#1}}
\newcommand{\bY}{\b{Y}}
\newcommand{\bYY}{\bY\mbox{$^*$}}
\newcommand{\bM}{\b{M}}
\newcommand{\mP}{\mbox{\textrm P}}

\title{New bounds for the Probability of Causation in Mediation
  Analysis} 

\author{Rossella Murtas\thanks{University of Cagliari} \and Alexander
  Philip Dawid\thanks{Leverhulme Emeritus Fellow, University of
    Cambridge} \and Monica Musio\thanks{University of Cagliari} }

\date{November 1, 2016}

\begin{document}
\maketitle

\begin{abstract}
  An individual has been subjected to some exposure and has developed
  some outcome.  Using data on similar individuals, we wish to
  evaluate, for this case, the probability that the outcome was in
  fact caused by the exposure.  Even with the best possible
  experimental data on exposure and outcome, we typically can not
  identify this ``probability of causation'' exactly, but we can
  provide information in the form of bounds for it.  Under appropriate
  assumptions, these bounds can be tightened if we can make other
  observations (e.g., on non-experimental cases), measure additional
  variables (e.g., covariates) or measure complete mediators. In this
  work we propose new bounds for the case that a third variable
  mediates partially the effect of the exposure on the outcome.\\

\end{abstract}

\section{Introduction}
\label{intro}
\label{sec:introduction}

Causality is a concept very common in real life situations.  Is lung
cancer caused by smoking?  Was contaminated water causing cholera in
London in 1854?  Can the court infer sex discrimination in a hiring
process?  However, statisticians have been very cautious in
formalizing this concept.  One reason may be the complex definitions
and methods implemented to study causality.  Another explanation may
be the difficulty of translating real life problems into mathematical
notations and formulas.  The first step should be to identify the
causal question of interest.  This can be assigned to one of two main
classes: questions concerning the causes of observed effects, and
questions concerning the effects of applied causes.  This basic
distinction, all too often neglected in the causal inference
literature, is fundamental to identifying the correct definition of
causation.  To clarify this distinction, consider the following
example.  An individual, Ann, might be subjected to some exposure,
$X$, and might develop some outcome, $Y$.  For simplicity we take $X$
to be a binary decision variable, denoting whether or not an
individual is given the drug, and take the outcome variable $Y$ also
to be binary, coded $1$ if the individual dies, and $0$ if not.  We
denote by $X_A\in\{0,1\}$ the value of Ann's exposure, and by
$Y_A\in\{0,1\}$ the value of Ann's outcome.  Questions about the
effects of applied causes, ``EoC'', are widely studied.  For example,
in medicine, randomized clinical trials are one of the most rigorous
ways to assess the effect of a treatment in a population.  In the EoC
framework, at an individual level we would be interested in asking:
``What would happen to Ann were she to be given the drug?'' or ``What
would happen to Ann were she not to be given the drug?''.  At the
population level, a typical EoC query would be: ``Is death caused by
the drug?''  In this framework, a straightforward way to assess the
strenght of causality is by comparing
$P_1=\mP(Y=1 \mid X \leftarrow 1)$ and
$P_0=\mP(Y=1 \mid X\leftarrow 0)$, the two outcome probabilities under
the two different interventions \cite{dawid2016statistical}.  This can
be seen as a decision problem: we can compare these two different
distributions for $Y$, decide which one we prefer, and take the
associated decision (give or withhold the drug).  The difference
$P_1 - P_0$ is known as the ``Average Causal Effect''.

In contrast to EoC queries, that are mostly adopted to infer knowledge
in the population, CoE questions invariably require an individual
investigation.  For example, suppose that Ann died after being given
the drug.  A typical CoE question might be phrased as: ``Knowing that
Ann did take the drug, and died, how likely is it that she would not
have died if she had not had the drug?''.  In this paper we will embed
such causal queries in the counterfactual framework
\cite{rubin1974estimating}.  This is based on the idea that there
exist \textit{potential variables}.  If $X$ is the exposure and $Y$
the outcome, the potential variable $Y(x)$ is conceived as the value
of $Y$ that would arise if, actually or hypothetically, $X$ were to be
set to $x$ ($X \leftarrow x$).  We denote the pair $(Y(0), Y(1))$ by
$\bY$.  For an actual assignment $X \leftarrow x$, we observe
$Y = Y(x)$.  The potential variable $Y(x')$, with $x' \neq x$, is then
not observable, but is supposed to describe what would have happened
to the outcome $Y$, if, counterfactually, we had assigned the
different value $x'$ to the exposure $X$.  Note particularly that it
is never possible to observe fully the pair $\bY$.

The definition of a CoE causal effect is completely different from the
EoC definition.  It is typically framed in terms of the
\textit{probability of causation} (PC), also called
\textit{probability of necessity} \cite{pearl1999probabilities}.
Given that Ann took the drug and died, the probability of causation in
Ann's case is defined as:
\begin{equation}\label{def:pcA}
  \pc_A = \mbox{P}_A(Y_A(0) = 0 \mid X_A= 1, Y_A(1)=1)
\end{equation}
where $\mbox{P}_A$ denotes the probability distribution over
attributes of Ann.  For example, suppose that Ann's children filed a
criminal lawsuit against a pharmaceutical manufacturer claiming that
their drug was the {\em cause\/} of her death.  Using data on similar
individuals, we would wish to evaluate, for this case, the probability
that the outcome was in fact caused by the exposure.

In such a civil case, the required standard of proof is typically
``preponderance of the evidence,'' or `` balance of probabilities,''
meaning that the case would succeed if it can be shown that causation
is ``more probable than not,'' \ie, $\pc_A > 50\%$.  However,
simplistic or {\em ad hoc\/} definitions and rules are widely and
often wrongly applied in many courthouse.  Given the possibly serious
implications of the probability of causation, it is important to
studying methods capable of producing accurate information.

From a statistical point of view, definition \eqref{def:pcA} involves
the bivariate distribution of the two potential variables associated
with the same subject.  However, only one of these can ever be
observed, the other then becoming counterfactual.  For this reason,
$\mbox{PC}_A$ is generally not fully identifiable.  We can however
provide useful information as bounds between which $\pc_A$ must lie.
Under appropriate assumptions, these bounds can be tightened if we can
measure additional variables, \eg, covariates \cite{dawid2011role},
or---in the case that unobserved variables confound the
exposure-outcome relationship---gather data on other, nonexperimental,
cases (Tian and Pearl \cite{tian2000probabilities}.)

In this paper we propose a novel approach to bound the probability of
causation in mediation analysis.  Mediation aims to disentangle the
extent to which the effect of $X$ on $Y$ is mediated through other
pathways from the extent to which that effect is due to $X$ acting
directly on $Y$.  In \secref{simple} we revisit the basic framework
where we have information only on exposure and outcome.  In
\secref{mediation} we focus on two different mechanisms: complete and
partial mediation.  In the former, the exposure is supposed to act on
the outcome only through the mediator, \ie, no direct effect is
present.  In the latter, both direct and indirect effects are
considered.  In \secref{comparison} we compare the bounds obtained in
\secref{mediation} with those reviewed in \secref{simple}, and in
\secref{conclusion} we present our conclusions.

\section{Starting Point: Simple Analysis}
\label{sec:simple}

In this Section we discuss the simple situation in which we have
information, as in \tabref{aspirin1}, from a randomized experimental
study that tested the same drug taken by Ann.

\begin{table}[htbp]
  \centering
  \begin{tabular}{lccc}
    \hline\noalign{\smallskip}
    & Die & Live & Total\\
    \noalign{\smallskip}\hline\noalign{\smallskip}
    Exposed & 30 & 70 & 100 \\
    Unexposed & 12 & 88 & 100 \\
    \noalign{\smallskip}\hline
  \end{tabular}
  \caption{Deaths in individuals exposed and unexposed to the same drug taken by Ann.}
  \label{tab:aspirin1}       
\end{table}


\begin{eqnarray}
  \mbox{P}(Y = 1 \mid X\leftarrow 1) &=& 0.30  \label{eq:pop-rate1}\\
  \mbox{P}(Y = 1 \mid X\leftarrow 0) &=& 0.12. \label{eq:pop-rate0}
\end{eqnarray}
We see that, in the experimental population, individuals exposed to
the drug ($X\leftarrow1$) were more likely to die than those unexposed
($X\leftarrow 0$), by $18$ percentage points.  So can the court infer
that it was Ann's taking the drug that caused her death?  More
generally: Is it correct to use such experimental results, concerning
a population, to say something about a single individual?  This
``Group-to-individual'' (G2i) issue is discussed by Dawid
\cite{dawid2014fitting}.  The simple difference between
\eqref{eq:pop-rate1} and \eqref{eq:pop-rate0} is not sufficient to
infer causation for a single external individual.

To make progress we add a further assumption that the event of Ann's
exposure, $X_A$, is independent of her potential response pair
$\textbf{Y}_A$:

\begin{equation}
  X_A \cip \bY_A.    \label{eq:suff}
\end{equation}

Property \eqref{eq:suff} parallels the ``no-confounding'' property
$X_i \cip \textbf{Y}_i$ that holds for individuals $i$ in the
experimental study on account of randomization.  We further suppose
that Ann is exchangeable with the individuals in the experiment, \ie,
she could be considered as a subject in the experimental population.
On account of \eqref{eq:suff} and exchangeability, $\mbox{PC}_A$ in
\eqref{def:pcA} reduces to
$\mbox{PC}_A = \mbox{P}(Y(0) = 0 \mid Y(1)=1)$---but we can not fully
identify this from the data.  In fact we can never observe the joint
event $(Y(0)=0;Y(1)=1)$, since at least one of $Y(0)$ and $Y(1)$ must
be counterfactual.  In particular, we can never learn anything about
the dependence between $Y(0)$ and $Y(1)$.  However, even without
making any assumptions about this dependence, we can derive the
following inequalities (Dawid \emph{et
  al.}~\cite{dawid2016statistical}):
\begin{equation}
  \label{eq:generic}
  \max\left\{0,1 - \frac{1}{\mbox{RR}}\right\} \leq
  \mbox{PC}_A
  \leq
  \frac{\min\{\mbox{P}(Y = 0 \mid X\leftarrow 0),\mbox{P}(Y = 1 \mid X\leftarrow 1)\}}
  {\mbox{P}(Y=1 \mid X\leftarrow 1)},
\end{equation}
where
\begin{equation}
  \mbox{RR} = \frac{\mbox{P}(Y = 1 \mid X\leftarrow 1)}{\mbox{P}(Y=1 \mid X\leftarrow 0)}
\end{equation}
is the {\em experimental risk ratio\/} between exposed and unexposed.
These bounds can be estimated from the experimental data using the
population death rates in equations~\eqref{eq:pop-rate1} and
\eqref{eq:pop-rate0}.

In many cases of interest (such as \tabref{aspirin1}), we have
\begin{displaymath}
  \mbox{P}(Y= 1 \mid X \leftarrow 0) <  \mbox{P}(Y= 1 \mid X \leftarrow 1)  < \mbox{P}(Y= 0 \mid X \leftarrow 0).
\end{displaymath}
Then the lower bound in \eqref{eq:generic} will be non-trivial, while
the upper bound will be 1, so vacuous.  Since, in \tabref{aspirin1},
the exposed are $2.5$ times as likely to die as the unexposed
($\mbox{RR}= 30/12 = 2.5$), we have enough confidence to infer
causality in Ann's case, since $0.60 \leq \mbox{PC}_A \leq 1$.

\section{Bounds in Mediation}
\label{sec:mediation}

In this Section we bound the probability of causation for a case where
a third variable, $M$, is involved in the causal pathway between the
exposure $X$ and the outcome $Y$.  We first review the results of
Dawid \emph{et al.}~\cite{dawid-murtas} for the case of complete
mediation, where no direct effect is present between exposure and
outcome but all the effect is mediated by $M$.  In addition we derive
new bounds for $\mbox{PC}_A$ in mediation analysis when a partial
mediation mechanism is in operation.

\subsection{Complete mediation (Dawid \emph{et
    al.}~\cite{dawid-murtas})}
\label{sec:compl-mediation}

\begin{figure}[h]
  \centering
  \begin{pspicture}(0,-0.5)(3,0.5) \color{black} \rput(0,0){%
      \rnode{1}{\Large X}} \rput(1.5,0){%
      \rnode{2}{\Large M}} \rput(3,0){%
      \rnode{3}{\Large Y}}
    \ncline[linecolor=black,nodesepA=3pt,nodesepB=3pt]{->}{1}{2}
    \ncline[linecolor=black,nodesepA=3pt,nodesepB=3pt]{->}{2}{3}
  \end{pspicture}
  \caption{Graph representing a mediator $M$, responding to exposure
    $X$ and affecting response $Y$.  There is no direct effect,
    unmediated by $M$, of $X$ on $Y$.}
  \label{fig:compl-med}
\end{figure}
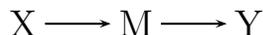

The case of no direct effect is intuitively described by
\figref{compl-med}.  Applications where this assumption might be
plausible is in the treatment of ovarian cancer (Silber \emph{et al.}
\cite{silber2007does}), where $X$ represents management either by a
medical oncologist or by a gynaecological oncologist, $M$ is the
intensity of chemotherapy prescribed, and $Y$ is death within 5 years.
We shall be interested in the case that $M$ is observed in the
experimental data but is not observed for Ann, and see how this
additional experimental evidence can be used to refine the bounds on
$\mbox{PC}_A$.

We now introduce $M(x)$ to denote the potential value of $M$ when
$X\leftarrow x$, and $Y^*(m)$ to denote the potential value of $Y$
when $M\leftarrow m$.  Then $Y(x):=Y^*\{M(x)\}$.  We define
$\bM: = (M(0),M(1))$ and $\bYY := (Y^*(0),Y^*(1))$.

We suppose that none of the causal mechanisms depicted in
\figref{compl-med} are confounded--expressed mathematically by
assuming mutual independence between $X$, $\bM$ and $\bYY$ (both for
experimental individuals, and for Ann). These assumptions imply no
overall confounding (as in \eqref{eq:suff}), the Markov property
$Y \cip X \mid M$, and the following bounds in this case of complete
mediation:

\begin{equation} \label{eq:pc-compl-med}
  \max\left\{0,{1-\frac{1}{\mbox{RR}}}\right\} \leq \mbox{PC}_A \leq
  \frac{\rm Num}{\mbox{P}(Y=1 \mid X\leftarrow 1)},
\end{equation}
where the numerator, Num, is given in \tabref{compl-med-upper}.  We
see from \eqref{eq:pc-compl-med} that knowing a mediator does not
improve the lower bound.  For the upper bound, one has to consider
various scenarios according to different choices of the estimable
marginal probabilities in \tabref{compl-med-upper}.

\begin{table}[htbp]
  \centering
  \begin{tabular}{l|cc}
    \hline\noalign{\smallskip}
    & $a\leq b$ & $a>b$\\
    \noalign{\smallskip}\hline\noalign{\smallskip}
    $c\leq d$ &  $a\cdot c+(1-d)(1-b)$ & $b\cdot c+(1-d)(1-a)$\\
    $c>d$ &  $a\cdot d+(1-c)(1-b)$ & $b\cdot d+(1-a)(1-c)$ \\
    \noalign{\smallskip}\hline
  \end{tabular}
  \caption{Numerator of upper bound for $\mbox{PC}_A$ in complete mediation anlaysis.  Here $a=\mP(M(0)=0)$, $b=\mP(M(1)=1)$, $c=\mP(Y^*(0)=0)$ and $d=P(Y^*(1)=1)$.}
  \label{tab:compl-med-upper}       
\end{table}

Note that, given the no-confounding assumptions, the entries in
\tabref{compl-med-upper} are all estimable from the experimental data:
\begin{eqnarray*}
  a &=& \mP(M = 0 \mid X \leftarrow 0)\\
  b &=& \mP(M = 1 \mid X \leftarrow 1)\\
  c &=& \mP(Y = 0 \mid M \leftarrow 0)\\
  d &=& \mP(Y = 1 \mid M \leftarrow 1).
\end{eqnarray*}

\subsection{Partial mediation}
\label{sec:part-mediation}
\begin{figure}[h]
  \centering
  \begin{pspicture}(-0.5,-0.5)(3.5,2) \color{black} \rput(0,0){%
      \rnode{1}{\LARGE X}} \rput(3,0){%
      \rnode{2}{\LARGE Y}} \rput(1.5,1.5){%
      \rnode{3}{\LARGE M}}
    \ncline[linecolor=black,nodesepA=3pt,nodesepB=3pt]{->}{1}{2}
    \ncline[linecolor=black,nodesepA=3pt,nodesepB=3pt]{->}{1}{3}
    \ncline[linecolor=black,nodesepA=3pt,nodesepB=3pt]{->}{3}{2}
  \end{pspicture}
  \caption{Graph illustrating a partial mediation mechanism between an
    exposure $X$, an outcome $Y$, and a mediator $M$}
  \label{fig:mediation-simple}
\end{figure}
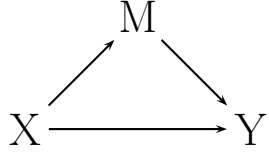


The situation described by \figref{compl-med} is unlikely to hold in
many real life situations.  Situations such as that represented
informally by \figref{mediation-simple}, that allow both direct and
indirect effects, are more plausible.  In this Section we derive new
bounds for the probability of causation when a partial mediator is
involved in the causal pathway.  We now define $Y^*(x,m)$ as the
potential value of the outcome $Y$ after setting both exposure,
$X\leftarrow x$, and mediator, $M \leftarrow m$.  Then
$Y(x)=Y^*(x,M(x))$.  We make the following assumptions, both for Ann
and for the individuals in the experiment:

\begin{description}
\item[A1:] $Y^*(x,m) \cip \bM \mid X$ (no $M$--$Y$ confounding)
\item[A2:] $Y^*(x,m) \cip X$ (no $X$--$Y$ confounding)
\item[A3:] $M(x) \cip X$ (no $X$--$M$ confounding).
\end{description}
Assumption \textbf{A1} expresses independence, given $X$, between a
potential value $Y^*(x,m)$, that would arise on setting exposure and
mediator to particular values, and the pair of potential outcomes
$(M(0), M(1))$.  It can be seen as a strengthening of the univariate
hypothesis $Y^*(x,m) \cip M(x) \mid X$.

Note that \textbf{A1} and \textbf{A2} are together equivalent to the
single requirement:
\begin{description}
\item[A12:] $Y^*(x,m) \cip (\bM, X)$.
\end{description}

Because we have supposed that Ann is exchangeable with the individuals
in the experiment, we have
\begin{equation}
  \label{eq:pceq}
  \mbox{PC}_A = \mbox{P}(Y(0) = 0 \mid X\leftarrow 1, Y(1)=1)=\frac{\mbox{P}(Y(0)=0,\: Y(1)=1 \mid X\leftarrow 1)}{\mbox{P}(Y(1)=1 \mid X\leftarrow 1)}.
\end{equation}

Given the no-confounding assumptions, the denominator of
\eqref{eq:pceq} is $\mP(Y = 1 \mid X \leftarrow 1)$, which is
estimable.  However, the numerator of \eqref{eq:pceq} involves the
joint distribution of the pair $\bY$ of potential outcomes, and this
is not estimable from the data, in view of the fact that it is never
possible to observe both $Y(0)$ and $Y(1)$ simultaneously.  We can
however bound this numerator in terms of estimable quantities, using
the fact that, for any events $A$ and $B$, and any probability
distribution $\mP$,
\begin{equation}\label{eq:bound-copula}
  \max\{\mP(A)+\mP(B)-1, 0 \} \leq \mP(A \cap B) \leq \min \{\mP(A),\mP(B)\}.
\end{equation}



Using \eqref{eq:bound-copula}, we can obtain an upper bound for the
numerator as: \scriptsize
\begin{align}
  &\mbox{P}(Y(0)=0,Y(1)=1 \mid X\leftarrow 1)=\mbox{P}(Y^*(0,M(0))=0,\:Y^*(1,M(1))=1 \mid X\leftarrow 1) \nonumber \\
    &=\sum_{m_0} \sum_{m_1} \mbox{P}(Y^*(0,m_0)=0, Y^*(1,m_1)=1 \mid M(0)=m_0,M(1)=m_1,X\leftarrow 1)
      \label{eq:8a}\\
  &\hspace{3cm} \times  \mbox{P}(M(0)=m_0,M(1)=m_1 \mid X\leftarrow 1) \nonumber \\
  &\leq \sum_{m_0} \sum_{m_1} \min\{\mbox{P}(Y^*(0,m_0)=0 \mid M(0)=m_0,M(1)=m_1,X\leftarrow 1),\nonumber\\
  &\hspace{3cm} \mbox{P}(Y^*(1,m_1)=1 \mid M(0)=m_0,M(1)=m_1,X\leftarrow 1)\} \nonumber \\
  & \hspace{3cm}\times \min\{\mbox{P}(M(0)=m_0 \mid X\leftarrow 1),\mbox{P}(M(1)=m_1 \mid X\leftarrow 1)\} \nonumber\\
  &=\sum_{m_0} \sum_{m_1}  \min\{\mbox{P}(Y^*(0,m_0)=0),\mbox{P}(Y^*(1,m_1)=1)\} \nonumber \\
  & \hspace{3cm}\times \min\{\mbox{P}(M(0)=m_0),\mbox{P}(M(1)=m_1\}, \label{eq:upper-part}
\end{align}
\normalsize on using assumptions \textbf{A12} and \textbf{A3}.  That
is, \footnotesize
\begin{align}
  &\mbox{P}(Y(0)=0,\:Y(1)=1 \mid X\leftarrow 1)\leq \nonumber \\
  &\min\{\mbox{P}(Y^*(0,0)=0),\mbox{P}(Y^*(1,0)=1)\} \times \min\{\mbox{P}(M(0)=0),\mbox{P}(M(1)=0)\} \label{eq:upper-part1}  \\
  &+\min\{\mbox{P}(Y^*(0,0)=0),\mbox{P}(Y^*(1,1)=1)\} \times \min\{\mbox{P}(M(0)=0),\mbox{P}(M(1)=1)\} \label{eq:upper-part2}\\
  &+\min\{\mbox{P}(Y^*(0,1)=0),\mbox{P}(Y^*(1,0)=1)\} \times \min\{\mbox{P}(M(0)=1),\mbox{P}(M(1)=0)\} \label{eq:upper-part3}\\
  &+\min\{\mbox{P}(Y^*(0,1)=0),\mbox{P}(Y^*(1,1)=1)\} \times \min\{\mbox{P}(M(0)=1),\mbox{P}(M(1)=1)\} \label{eq:upper-part4}.
\end{align}
\normalsize

It can be shown that similarly applying the lower bound of
\eqref{eq:bound-copula} to \eqref{eq:8a} yields the same lower bound
as obtained in \secref{simple} and \secref{compl-mediation}, so that
the lower bound is not improved by knowledge of a mediation mechanism.

Assumptions \textbf{A12} and \textbf{A3} allow us to estimate the
terms \eqref{eq:upper-part1}--\eqref{eq:upper-part4} in the above
upper bound from the data:
\begin{eqnarray*}
  \mP(Y^*(x,m)=y) &=& \mP(Y=y \mid X\leftarrow x, M\leftarrow m)\\
  \mP(M(x) = m) &=& \mP(M = m \mid X\leftarrow x).
\end{eqnarray*}





\section{Comparisons}
\label{sec:comparison}

In this Section we compare the bounds found in the simple analysis
framework of \secref{simple} with those obtained by considering a
complete mediation mechanism, as in \secref{compl-mediation}, and
those obtained by considering a partial mediation mechanism, as in
\secref{mediation}.  We focus on comparing these bounds to obtain the
best information from the data.

The numerator of the upper bound for $\mbox{PC}_A$ in
\eqref{eq:generic}, which ignores the mediator, may be written as
\begin{equation}
  \label{eq:alpha}
  \min\{\alpha+\beta,\gamma+\delta\},
\end{equation}
where
\begin{eqnarray*}
  \alpha &=& \mP(Y^*(0,0)=0)\,\mP(M(0)=0)\\
  \beta &=& \mP(Y^*(0,1)=0)\,\mP(M(0)=1)\\
  \gamma &=& \mP(Y^*(1,0)=1)\,\mP(M(1)=0)\\
  \delta &=& \mP(Y^*(1,1)=1)\,\mP(M(1)=1).
\end{eqnarray*}

We see that both \eqref{eq:upper-part1} and \eqref{eq:upper-part2} are
smaller than or equal to $\alpha$, while both \eqref{eq:upper-part3}
and \eqref{eq:upper-part4} are smaller than or equal to $\beta$.  So
the upper bound allowing for partial mediation, which is the sum of
\eqref{eq:upper-part1}, \eqref{eq:upper-part2}, \eqref{eq:upper-part3}
and \eqref{eq:upper-part4}, cannot exceed
$2\,(\alpha+\beta) = 2\,\mP(Y(0)=0) = 2\, \mP(Y=0 \mid X \leftarrow
0)$; and similarly cannot exceed $2\,\mP(Y=1 \mid X \leftarrow 1)$.
Thus, the upper bound for the numerator, when accounting for the
mediator, can not exceed twice that obtained by ignoring it, as given
by \eqref{eq:generic}.  However, as we will see in \secref{example},
it could be larger or smaller than that simpler bound.  On the other
hand, we do not obtain a different lower bound.


In the special case of complete mediation, $Y^*(0,m) = Y^*(1,m)$,
$= Y^*(m)$, say.  Thus the terms with $m_0\neq m_1$ in \eqref{eq:8a}
must be $0$.  This leads to the following upper bound:
\begin{align}
  &\mbox{P}(Y(0)=0,Y(1)=1 \mid X\leftarrow 1)\leq  \nonumber \\
  &\min\{\mbox{P}(Y^*(0)=0),\mbox{P}(Y^*(1)=1)\} \times \min\{\mbox{P}(M(0)=0),\mbox{P}(M(1)=1)\}+ \nonumber \\
  &+\min\{\mbox{P}(Y^*(1)=0),\mbox{P}(Y^*(0)=1)\} \times \min\{\mbox{P}(M(0)=1),\mbox{P}(M(1)=0)\}, \nonumber
\end{align}
in agreement with \tabref{compl-med-upper}.  Since we have eliminated
the terms \eqref{eq:upper-part1} and \eqref{eq:upper-part4} appearing
in the general case, the upper bound obtained in this case of complete
mediation is never bigger than that obtained from the general
expression
\eqref{eq:upper-part1}+\eqref{eq:upper-part2}+\eqref{eq:upper-part3}+\eqref{eq:upper-part4};
nor, since \eqref{eq:upper-part2} $\leq \alpha$ while
\eqref{eq:upper-part3} $\leq \beta$, can it be bigger than the bound
\eqref{eq:generic} obtained on ignoring the knowledge of the complete
mediation mechanism.

\subsection{Examples}
\label{sec:example}

To show that, in the case of partial mediation, taking account of
information about the mediator may, but need not, yield a tighter
upper bound, we consider two cases with different experimental data as
given respectively by \tabref{aspirin2} and \tabref{aspirin3}.

\begin{table}[htbp]
  \centering
  \begin{tabular}{lccc}
    \hline\noalign{\smallskip}
    & Die & Live & Total\\
    \noalign{\smallskip}\hline\noalign{\smallskip}
    Exposed & 69 & 31 & 100\\
    Unexposed & 24 & 76 & 100\\
    \noalign{\smallskip}\hline
  \end{tabular}
  \caption{Experimental data~1}
  \label{tab:aspirin2}       
\end{table}

\begin{table}[htbp]
  \centering
  \begin{tabular}{lccc}
    \hline\noalign{\smallskip}
    & Die & Live & Total\\
    \noalign{\smallskip}\hline\noalign{\smallskip}
    Exposed & 78 & 22 & 100\\
    Unexposed & 32 & 68 & 100\\
    \noalign{\smallskip}\hline
  \end{tabular}
  \caption{Experimental data~2}
  \label{tab:aspirin3}       
\end{table}

Suppose now we can also observe a partial mediator $M$.  We might then
observe the following probabilities, consistent with
\tabref{aspirin2}:

\begin{align}
  & \mbox{P}(Y^*(0,0)=0)= 0.98  & \mbox{P}(Y^*(0,1)=0) = 0.165 \nonumber \\
  & \mbox{P}(Y^*(1,0)=0)= 0.315 & \mbox{P}(Y^*(1,1)=0) = 0.143 \nonumber\\
  & \mbox{P}(M(0)=0) = 0.73   & \mbox{P}(M(1)=0) = 0.981 \nonumber
\end{align}

We then obtain: $0.65 \leq \mbox{PC}_A \leq 0.81$ when accounting for
the mediator, and $0.65 \leq \mbox{PC}_A \leq 1$ when ignoring it.  In
this case, knowledge of the partial mediation mechanism is helpful in
improving the upper bound.

On the other hand, suppose we observe the following probabilities,
consistent with \tabref{aspirin3}:

\begin{align}
  & \mbox{P}(Y^*(0,0)=0)= 0.98  & \mbox{P}(Y^*(0,1)=0) = 0.67 \nonumber \\
  & \mbox{P}(Y^*(1,0)=0)= 0.09 & \mbox{P}(Y^*(1,1)=0) = 0.27 \nonumber\\
  & \mbox{P}(M(0)=0) = 0.04   & \mbox{P}(M(1)=0) = 0.26 \nonumber
\end{align}

We now obtain: $0.59 \leq \mbox{PC}_A \leq 0.95 $ when accounting for
the mediator, but $0.59 \leq \mbox{PC}_A \leq 0.88$ when ignoring it.
So in this case knowledge of mediation has not been helpful.

\section{Conclusions}
\label{sec:conclusion}

Bounding the probability of causation in mediation analysis is an
important problem for applications.  By taking account of a complete
mediation mechanism we can never do worse than by ignoring it.
However, complete mediation is not always reasonable.  In the case of
partial mediation, the upper bound obtained by taking account of it
may be greater or smaller than that obtained by ignoring it.  We can
thus compute both upper bounds and take the smaller.

This work has several possible extensions.  It would be interesting to
extend our theoretical formulas to cases combining information on both
covariates and mediators.  Another promising extension arises on
making connections with copula theory, where $\mbox{PC}$ can be
obtained as a function of the estimable quantities $\mP(Y(0)=0)$ and
$\mP(Y(1)=1)$ together with an appropriate copula.


\begin{thebibliography}{}


\bibitem{dawid2011role} Dawid AP (2011) The role of scientific and
  statistical evidence in assessing causality. In: Goldberg~R (ed)
  Perspectives on Causation. Oxford, Hart Publishing, pp 133--147

\bibitem{dawid2014fitting} Dawid AP, Faigman DL, Fienberg SE (2014)
  Fitting science into legal contexts: Assessing effects of causes or
  causes of effects?\@ Sociological Methods \& Research 43:359--390

\bibitem{dawid-murtas} Dawid AP, Murtas R, Musio M (2016) Bounding the
  probability of causation in mediation analysis. In: Di Battista T,
  Moreno E, Racugno W (ed) Topics on Methodological and Applied
  Statistical Inference.  Cham Switzerland, Springer International
  Publishing, pp 75--84

\bibitem{dawid2016statistical} Dawid AP, Musio M, Fienberg SE (2016)
  From statistical evidence to evidence of causality. Bayesian
  Analysis, 11(3):725--752
%

\bibitem{silber2007does} Markman M (2007) Does ovarian cancer
  treatment and survival differ by the specialty providing
  chemotherapy? Journal of Clinical Oncology, 25(23):3554--3554

\bibitem{pearl1999probabilities} Pearl J (1999) Probabilities of
  causation: Three counterfactual interpretations and their
  identification. Synth\`ese, 121(1--2):93--149


\bibitem{rubin1974estimating} Rubin DB (1974) Estimating causal
  effects of treatments in randomized and nonrandomized
  studies. Journal of Educational Psychology, 66(5):688

\bibitem{tian2000probabilities} Tian J, Pearl J (2000) Probabilities
  of causation: Bounds and identification. Annals of Mathematics and
  Artificial Intelligence, 28(1--4):287--313

\end{thebibliography}
\end{document}